\documentclass[12pt,twoside]{article}
\usepackage{graphicx}
\usepackage{amsfonts,amsbsy,amssymb,amsmath}
\allowdisplaybreaks
\usepackage{cite}
\usepackage{graphics}
\def\bpsp{\begin{pspicture}}
\def\epsp{\end{pspicture}}

\newtheorem{theorem}{Theorem}[section]
\newtheorem{remark}[theorem]{Remark}
\newtheorem{example}[theorem]{Example}
\newtheorem{lemma}[theorem]{Lemma}
\newtheorem{corollary}[theorem]{Corollary}
\newtheorem{definition}[theorem]{Definition}
\newtheorem{proposition}[theorem]{Proposition}

\newtheorem{note}{Note}
\newtheorem{case}{Case}
\newtheorem{conjecture}{Conjecture}
\newtheorem{question}{Question}

\newcommand{\bea}{\begin{eqnarray}}
\newcommand{\eea}{\end{eqnarray}}
\newcommand{\beq}{\begin{eqnarray*}}
\newcommand{\eeq}{\end{eqnarray*}}

\def\m4{\mbox{\rm ~(mod $4$)}}

\def \bd{\begin{definition}}
\def \ed{\end{definition}}
\def \bqu{\begin{question}}
\def \equ{\end{question}}
\def \bcc{\begin{conjecture}}
\def \ecc{\end{conjecture}}
\def \bt{\begin{theorem}}
\def \et{\end{theorem}}
\def \bl{\begin{lemma}}
\def \el{\end{lemma}}
\def \bc{\begin{corollary}}
\def \ec{\end{corollary}}
\def \be{\begin{equation}}
\def \ee{\end{equation}}
\def \ben{\begin{enumerate}}
\def \een{\end{enumerate}}
\def \ba{\begin{array}}
\def \ea{\end{array}}
\def \bp{\begin{proposition}}
\def \ep{\end{proposition}}
\def \bx{\begin{example}}
\def \ex{\end{example}}
\def \br{\begin{remark}}
\def \er{\end{remark}}
\def \bdsc{\begin{description}}
\def \edsc{\end{description}}

\def \bn{\begin{case}}
\def \en{\end{case}}
\def \bnt{\begin{note}}
\def \ent{\end{note}}
\def\1{1\!\!1}

\def\mm2{\mbox{\rm ~(mod $2$)}}
\def\m4{\mbox{\rm ~(mod $4$)}}

\def\qed{\nolinebreak\hfill\rule{.2cm}{.2cm}\par\addvspace{.5cm}}

\def\m{\mu}

\def\1{\textbf{1}}
\def\0{\textbf{0}}

\begin{document}
\title{On Zagreb index, signless Laplacian eigenvalues and signless Laplacian energy of a graph}
\author{S. Pirzada$ ^{a} $, Saleem Khan$ ^{b} $ \\
$^{a,b}${\em Department of Mathematics, University of Kashmir, Srinagar, Kashmir, India}\\
$^{a}$\texttt{pirzadasd@kashmiruniversity.ac.in}; $ ^{b} $\texttt{khansaleem1727@gmail.com}
}
\date{}

\pagestyle{myheadings} \markboth{Pirzada,Saleem}{On Zagreb index, signless Laplacian eigenvalues and energy of a graph}
\maketitle
\vskip 5mm
\noindent{\footnotesize \bf Abstract.} Let $G$ be a simple graph with order $n$ and size $m$. The quantity $M_1(G)=\displaystyle\sum_{i=1}^{n}d^2_{v_i}$ is called the first Zagreb index of $G$, where $d_{v_i}$ is the degree of vertex $v_i$, for all $i=1,2,\dots,n$. The signless Laplacian matrix of a graph $G$ is $Q(G)=D(G)+A(G)$, where $A(G)$ and $D(G)$ denote, respectively, the adjacency and the diagonal matrix of the vertex degrees of $G$. Let $q_1\geq q_2\geq \dots\geq q_n\geq 0$ be the signless Laplacian eigenvalues of $G$. The largest signless Laplacian eigenvalue $q_1$ is called the signless Laplacian spectral radius or $Q$-index of $G$ and is denoted by $q(G)$.  Let $S^+_k(G)=\displaystyle\sum_{i=1}^{k}q_i$ and $L_k(G)=\displaystyle\sum_{i=0}^{k-1}q_{n-i}$, where $1\leq k\leq n$, respectively denote the sum of $k$ largest and smallest signless Laplacian eigenvalues of $G$. The signless Laplacian energy of $G$ is defined as $QE(G)=\displaystyle\sum_{i=1}^{n}|q_i-\overline{d}|$, where $\overline{d}=\frac{2m}{n}$ is the average vertex degree of $G$. In this article, we obtain upper bounds for the first Zagreb index $M_1(G)$ and show that each bound is best possible. Using these bounds, we obtain several upper bounds for the graph invariant $S^+_k(G)$ and characterize the extremal cases. As a consequence,  we find upper bounds for the $Q$-index and lower bounds for the graph invariant $L_k(G)$ in terms of various graph parameters and determine the extremal cases. As an application, we obtain upper bounds for the signless Laplacian energy of a graph and characterize the extremal cases.

\vskip 3mm



\section{Introduction} We consider simple graphs $G=G(V,E)$ with order $n$ and size $m$ having vertex set $V(G)=\{v_1,v_2,\dots,v_n\}$ and edge set $E(G)=\{e_1,e_2,\dots,e_m\}$. As usual, $K_{1,n-1}$ and $K_n$ denote the star on $n$ vertices and the complete graph on $n$ vertices, respectively. The degree of a vertex $v_i\in V(G)$, denoted by $d_{v_i}=d_i$, is the number of edges incident on $v_i$. We will denote by $\triangle(G)$ and $\delta(G)$ the maximum vertex  degree and the minimum vertex degree in a graph $G$, respectively. The diameter of a connected graph $G$, denoted by $D(G)$, is the largest distance between any pair of vertices in $G$.  We refer the reader to \cite{R22,R23} for other undefined notations and terminology from spectral graph theory.

The adjacency matrix $A(G)=(a_{ij})$ of $G$ is a $(0,1)$-square matrix of order $n$ whose $(i,j)$-entry is equal to 1, if $v_i$ is adjacent to $v_j$ and equal to $0$, otherwise. If $\lambda_1\geq \lambda_2\geq\dots \geq \lambda_n$ are the adjacency eigenvalues of $G$, the energy \cite{R6} of $G$ is defined as $E(G)=\displaystyle\sum_{i=1}^{n}|\lambda_i|$. The quantity $E(G)$ introduced by I. Gutman has well developed mathematical aspect and has noteworthy chemical applications (see \cite{R7}).

Let $D(G)=diag(d_1,d_2,\dots,d_n)$ be the diagonal matrix associated to $G$, where $d_i=d_{v_i}$ is the degree of the vertex $v_i$, for all $i=1,2,\dots,n$. The matrices $L(G)=D(G)-A(G)$ and $Q(G)=D(G)+A(G)$ are called the Laplacian and the signless Laplacian matrices, respectively. Their spectrum are called the Laplacian spectrum and the signless Laplacian spectrum of the graph $G$, respectively. Both the matrices $L(G)$ and $Q(G)$ are real symmetric, positive semi-definite matrices, therefore their eigenvalues are non-negative real numbers. Let $\mu_1\geq \mu_2\geq \dots \geq \mu_n\geq 0$ and $ q_1\geq q_2\geq \dots \geq q_n\geq 0$ be the Laplacian spectrum and the signless Laplacian spectrum of the graph $G$, respectively. The eigenvalues of $Q(G)$ are called the $Q$-eigenvalues of $G$. Also, the largest signless Laplacian eigenvalue $q_1$ of $Q(G)$ is called the signless Laplacian spectral radius or $Q$-index of $G$ and is denoted by $q(G)$. For $k=1,2,\dots,n$, let
$S_k(G)=\displaystyle\sum_{i=1}^{k}\mu_i,$ be the sum of $k$ largest Laplacian eigenvalues of $G$. We note that the sum $S_k(G)$ is of much interest by itself and some exciting details, extensions and open problems about it may be found in the excellent paper of Nikiforov \cite{R8}. The well-known Brouwer's conjecture, due to Brouwer \cite{R9} about the sum $S_k(G)$ is stated as follows.
\begin{conjecture} If $G$ is any graph with order $n$ and size $m$, then
\begin{equation*}
S_k(G)\leq m+{k+1\choose 2},~~~~~~ \text{for any} ~~~ k\in \{1,2,\dots,n\}.
\end{equation*}
\end{conjecture}

Although Conjecture 1  has been studied extensively but it remains open at large. For the progress on Brouwer's Conjecture, we refer to \cite{R10,R11,R13} and the references therein.

Let $S^+_k(G)=\displaystyle\sum_{i=1}^{k}q_i$ and $L_k(G)=\displaystyle\sum_{i=0}^{k-1}q_{n-i}$, where $k=1,2,\dots,n$, be the sum of $k$ largest and smallest signless Laplacian eigenvalues of $G$, respectively.
Motivated by the studies of Mohar \cite{R14}, Jin et al. \cite{R15} investigated the sum of the $k$ largest signless Laplacian eigenvalues. Motivated by the definition of $S_k(G)$ and Brouwer's conjecture, Ashraf et al. \cite{R16} proposed the following conjecture about $S^+_k(G)$.
\begin{conjecture} If $G$ is any graph with order $n$ and size $m$, then
\begin{equation*}
S^+_k(G)\leq m+{k+1\choose 2},~~~~~~ \text{for any} ~~~ k\in \{1,2,\dots,n\}.
\end{equation*}
\end{conjecture}
To see the progress on this conjecture, we refer to \cite{R17} and the references therein.

The rest of the paper is organized as follows. In Section 2,  we obtain upper bounds for the first Zagreb index $M_1(G)$ and show that the bounds are sharp. Using these investigations, we obtain several upper bounds for the graph invariant $S^+_k(G)$ and determine the extremal graphs. As a consequence, we obtain upper bounds for the $Q$-index and lower bounds for the graph invariant $L_k(G)$ in terms of various graph parameters and determine the extremal cases in each case. In Section 3, we find some upper bounds for the signless Laplacian energy $QE(G)$ for a connected graph $G$ and determine the extremal cases.

\section{Sum of the signless Laplacian eigenvalues of a graph}

The first Zagreb index $M_1 (G)$ \cite{R24} of a graph $G$ is defined as
$M_1(G)=\displaystyle\sum_{i=1}^{n}d^2_{v_i}$, where $d_{v_i}$ is the degree of vertex $v_i$, for all $i=1,2,\dots,n.$ The following inequality can be found in \cite{R1}.

\begin{lemma} \label{L1}\emph{\textbf{\cite{R1}}} Let $a=(a_1 ,a_2 \dots,a_n)$ and $b=(b_1 ,b_2 \dots,b_n)$ be $n$-tuples of real numbers satisfying $0\leq m_1 \leq a_i \leq M_1$, $0\leq m_2\leq b_i \leq M_2$ with $i=1,2,\dots,n$ and $M_1 M_2 \neq=0$. Let $\alpha=\frac{m_1}{M_1}$ and $\beta=\frac{m_2}{M_2}$. If $(1+\alpha)(1+\beta)\geq 2$, then \begin{equation}\label{E1} \displaystyle\sum_{i=1}^{n}a_i^2 \displaystyle\sum_{i=1}^{n}b_i^2-{\Big(\displaystyle\sum_{i=1}^{n}a_i b_i\Big)}^2 \leq \frac{n^2}{4}{(M_1 M_2 -m_1 m_2)}^2. \end{equation}
\end{lemma}

The following result gives an upper bound for the graph invariant $M_1(G)$  in terms of the order $n$, size $m$, $\triangle(G)$ and $\delta(G)$.

\begin{lemma} \label{L5}
Let $G$ be a connected graph with $n$ vertices and $m$ edges. Then
\begin{equation}\label{E2}
\displaystyle\sum_{i=1}^{n}d_i^2\leq \frac{4m^2}{n}+\frac{n}{4}{(\triangle(G) -\delta(G))}^2.
\end{equation}
Furthermore, the inequality is sharp and is shown by all degree regular graphs.
\end{lemma}
\textbf{Proof.} In Lemma \ref{L1}, taking $a=(d_1,d_2,\dots,d_n)$, $b=(1,1,\dots,1)$, $M_1=\triangle(G)$, $m_1=\delta(G)$ and $M_2=m_2=1$. With these values the condition $(1+\alpha)(1+\beta)\geq 2$ in Lemma \ref{L1} is satisfied. Substituting these values in Inequality \ref{E1}, we get
\begin{equation*}
\displaystyle\sum_{i=1}^{n}d_i^2 \displaystyle\sum_{i=1}^{n}1-{\Big(\displaystyle\sum_{i=1}^{n}d_i\Big)}^2 \leq \frac{n^2}{4}{(\triangle(G) -\delta(G))}^2.
\end{equation*}
Using the fact that $\displaystyle\sum_{i=1}^{n}d_i=2m$ in the above inequality and simplifying further, we get
$$n\displaystyle\sum_{i=1}^{n}d_i^2-4m^2\leq \frac{n^2}{4}{(\triangle(G) -\delta(G))}^2,$$
that is,
$$\displaystyle\sum_{i=1}^{n}d_i^2\leq \frac{4m^2}{n}+\frac{n}{4}{(\triangle(G) -\delta(G))}^2,$$
which proves the required inequality.

Now, let $G$ be an $r$-regular graph so that $\triangle(G) =\delta(G)$. Clearly, the left hand side of Inequality \ref{E2} becomes $nk^2$ and the right hand side becomes $\frac{n^2k^2}{n}=nk^2$. This completes the proof. \qed

The next lemma shows that the diameter of a connected graph $G$ can be at most $e(G)-1$ where $e(G)$ is the number of distinct $Q$-eigenvalues of $G$.
\begin{lemma} \label{L4} Let $G$ be a connected graph of diameter $D$ and $e(G)$ distinct $Q$-eigenvalues. Then $D\leq e(G)-1$.
\end{lemma}

In the next lemma, we show that the complete graph is the unique connected graph having only two distinct $Q$-eigenvalues.

\begin{lemma} \label{L6}
 Let $G$ be a connected graph on $n$ vertices with $e(G)$ distinct $Q$-eigenvalues. Then $e(G)=2$ if and only if $G\cong K_n$.
 \end{lemma}
 \textbf{Proof.} Assume that $e(G)=2$. Then, from Lemma \ref{L4}, we have $D(G)=1$, which shows that $G\cong K_n$.

 Conversely, suppose that $G\cong K_n$. The proof follows by observing that the $Q$-spectrum of $K_n$ is $\{2n-2, \underbrace{n-2,\dots,n-2}_{n-1}\}$. \qed

A simpler version of classical Cauchy- Schwarz Inequality is as follows.

\begin{lemma} \label{L3} Let $(a_1,a_2,\dots,a_n)$ be a sequence of non-negative real numbers. Then $$\Big(\displaystyle\sum_{i=1}^{n}a_i\Big)^2\leq n\displaystyle\sum_{i=1}^{n}a^2_i$$ with equality if and only if $a_1=a_2=\dots=a_n$.
\end{lemma}

Now, we obtain an upper bound for $S^+_k(G)$ in terms of $n$, $m$, $\triangle(G)$ and $\delta(G)$ and characterize the extremal graphs.

 \begin{theorem}\label{T1}
 Let $G$ be a connected graph with $n$ vertices and $m$ edges. If $1\leq k\leq n-1$, then
 \begin{equation}\label{E3}
 S^+_k(G)\leq \frac{2mk}{n}+\frac{\sqrt{k(n-k)\Big(8mn+n^2{(\triangle(G)-\delta(G))}^2\Big)}}{2n}
 \end{equation}
 with equality if and only if $G\cong K_n$ and $k=1$. Equality always holds when $k=n$.
\end{theorem}
\textbf{Proof.} Using the fact that the sum of the eigenvalues of a matrix equals its trace, we have
$$2m=\displaystyle\sum_{v_i\in V(G)}d_{v_i}=q_1+q_2+\dots+q_n,$$
that is, $$2m+\displaystyle\sum_{v_i\in V(G)}d^2_{v_i}=\displaystyle\sum_{v_i\in V(G)}(d^2_{v_i}+d_{v_i})=q^2_1+q^2_2+\dots+q^2_n. $$
Let $S^+_k(G)=S^+_k$. Using the above equations with Lemma \ref{L3}, we get
\begin{align*}
(q_{k+1}+\dots+q_n)^2 &= (2m-S^+_k)^2 \leq (n-k)(q^2_{k+1}+\dots+q^2_n)\\&
= (n-k)\Big(2m+\displaystyle\sum_{v_i\in V(G)}d^2_{v_i}-(q^2_1+\dots+q^2_k)\Big)\\&
\leq  (n-k)\Big(2m+\displaystyle\sum_{v_i\in V(G)}d^2_{v_i}-\frac{{S^+_k}^2}{k}\Big).
\end{align*}
Simplifying further, we get
$${S^+_k}^2-\frac{4mkS^+_k}{n}+\frac{4m^2k}{n}-\frac{k(n-k)}{n}\Big(2m+\displaystyle\sum_{v_i\in V(G)}d^2_{v_i}\Big)\leq 0,$$
that is,
$$S^+_k\leq \frac{2mk}{n}+\frac{\sqrt{4m^2k^2-4knm^2+nk(n-k)\Big(2m+\displaystyle\sum_{v_i\in V(G)}d^2_{v_i}\Big)}}{n}$$
or
\begin{equation}\label{E4}
S^+_k\leq \frac{2mk}{n} +\frac{\sqrt{k(n-k)\Big(n(2m+\displaystyle\sum_{v_i\in V(G)}d^2_{v_i})-4m^2\Big)}}{n}.
\end{equation}
Using Lemma \ref{L5} in Inequality (\ref{E4}), we get
$$S^+_k\leq \frac{2mk}{n} +\frac{\sqrt{k(n-k)\Big(2mn+4m^2+\frac{n^2}{4}{(\triangle(G)-\delta(G))}^2-4m^2\Big)}}{n}$$
or
$$S^+_k\leq \frac{2mk}{n} +\frac{\sqrt{k(n-k)\Big(8mn+n^2{(\triangle(G)-\delta(G))}^2\Big)}}{2n}$$
and this proves the required inequality.

 Now, suppose that the equality holds in Inequality \ref{E3}. Then, from the above proof, equality must hold in Lemma \ref{L3} and Lemma \ref{L5}. Thus, we must have $q_{k+1}=q_{k+2}=\dots=q_n$ and $q_1=q_2=\dots=q_k$, from Lemma \ref{L3}. These two equalities show that $G$ has exactly two distinct $Q$-eigenvalues. Thus, by Lemma \ref{L6}, $G\cong K_n$ and we know that $K_n$ is a regular graph. Lastly, $k=1$ follows from the $Q$-spectrum of $K_n$.

 Conversely, it is easy to see that the equality holds in Inequality \ref{E3} if $G\cong K_n$ and $k=1$.

 Furthermore, if $k=n$ then the left hand side of Inequality \ref{E3} is $q_1+\dots+q_n=2m$ and the right hand side becomes $\frac{2mn}{n}=2m$. Thus, equality always holds when $k=n$. \qed

Proceeding and using arguments similar to those used in Theorem \ref{T1}, we get the following lower bound for $L_k(G)$.

\begin{theorem}\label{T2}
 Let $G$ be a connected graph with $n$ vertices and $m$ edges. If $1\leq k\leq n-1$, then
 \begin{equation*}
 L_k(G)\geq \frac{2mk}{n}-\frac{\sqrt{k(n-k)\Big(8mn+n^2{(\triangle(G)-\delta(G))}^2\Big)}}{2n}
 \end{equation*}
 with equality if and only if $G\cong K_n$ and $k=n-1$. Equality always holds when $k=n$.
\end{theorem}

Taking $k=1$ in Theorem \ref{T1}, we obtain the following upper bound for the signless Laplacian spectral radius in terms of $m$, $n$, $\triangle(G)$ and $\delta(G$.

\begin{theorem}\label{T3}
 Let $G$ be a connected graph with $n$ vertices and $m$ edges. Then
 \begin{equation*}
 q(G)\leq \frac{2m}{n}+\frac{\sqrt{(n-1)\Big(8mn+n^2{(\triangle(G)-\delta(G))}^2\Big)}}{2n}
 \end{equation*}
 with equality if and only if $G\cong K_n$.
\end{theorem}

The following inequality can be seen in \cite{R4}.
\begin{lemma}\label{L7} \emph{ \cite{R4}}
  If $a_i$ and $b_i$, $1\leq i\leq n$, are positive real numbers, then
  $$\displaystyle\sum_{i=1}^{n}a_i^2\displaystyle\sum_{i=1}^{n}b_i^2\leq \frac{1}{4}{\Bigg(\sqrt{\frac{M_1 M_2}{m_1 m_2}}+\sqrt{\frac{m_1 m_2}{M_1 M_2}}\Bigg)}^2\Bigg{(\displaystyle\sum_{i=1}^{n}a_i b_i\Bigg)}^2,$$
  where $M_1 =\text{max} \{a_i:1\leq i\leq n\}$,  $m_1 =\text{min} \{a_i:1\leq i\leq n\}$,  $M_2 =\text{max} \{b_i:1\leq i\leq n\}$ and  $m_2 =\text{min} \{b_i:1\leq i\leq n\}$.
  \end{lemma}

Now, we obtain a different upper bound for the sum of squares of the vertex degrees of a connected graph $G$ in terms of the same parameters as in Lemma \ref{L5}.

  \begin{lemma} \label{L8}
Let $G$ be a connected graph with $n$ vertices and $m$ edges. Then
\begin{equation}\label{E5}
\displaystyle\sum_{i=1}^{n}d_i^2\leq \frac{m^2\Big(\triangle(G) +\delta(G)\Big)^2}{n\triangle(G) \delta(G)}.
\end{equation}
Moreover, the inequality is sharp and is shown by all degree regular graphs.
\end{lemma}
\textbf{Proof.} In Lemma \ref{L7}, take $a_i=d_{v_i}=d_i$ ($1\leq i\leq n$), $b_i=1$ ($1\leq i\leq n$), $M_1=\triangle(G)$, $m_1=\delta(G)$ and $M_2=m_2=1$, we get
$$\displaystyle\sum_{i=1}^{n}d_i^2\displaystyle\sum_{i=1}^{n}1\leq \frac{1}{4}{\Bigg(\sqrt{\frac{\triangle(G)}{\delta(G)}}+\sqrt{\frac{\delta(G)}{\triangle(G)}}\Bigg)}^2{\Bigg(\displaystyle\sum_{i=1}^{n}d_i\Bigg)}^2. $$
Using $\displaystyle\sum_{i=1}^{n}d_i=2m$ in the above inequality, we get
$$\displaystyle\sum_{i=1}^{n}d_i^2\leq \frac{m^2\Big(\triangle(G) +\delta(G)\Big)^2}{n\triangle(G) \delta(G)},$$
which is the required inequality.

For the equality part, let $G$ be $k$-regular. Then the left hand side of Inequality \ref{E5} becomes $nk^2$ and the right hand side becomes $\frac{4k^4n^2}{4nk^2}=nk^2$. Thus equality holds in Inequality \ref{E5} whenever $G$ is a regular graph. \qed

Now, we will use the Lemma \ref{L8} to get the following upper bound for the graph invariant $S^+_k(G)$.
\begin{theorem}\label{T4}
Let $G$ be a connected graph with $n$ vertices and $m$ edges. If $1\leq k\leq n-1$, then
 \begin{equation}\label{E6}
 S^+_k(G)\leq \frac{2mk}{n}+\frac{\sqrt{mk(n-k)\Big(2\triangle(G)\delta(G)(n-2m)+m{(\triangle(G)+\delta(G))}^2\Big)}}{n\sqrt{\triangle(G)\delta(G)}}
 \end{equation}
 with equality if and only if $G\cong K_n$ and $k=1$. Equality always holds when $k=n$.
\end{theorem}
\textbf{Proof.} Proceeding similarly as in Theorem \ref{T1} upto Inequality \ref{E4} and using Lemma \ref{L8}, we get
$$S^+_k\leq \frac{2mk}{n} +\frac{\sqrt{k(n-k)\Big(n\Big(2m+\frac{m^2(\triangle(G) +\delta(G))^2}{n\triangle(G) \delta(G)}\Big)-4m^2\Big)}}{n}$$
or
$$S^+_k\leq \frac{2mk}{n}+\frac{\sqrt{mk(n-k)\Big(2\triangle(G) \delta(G)(n-2m)+m{(\triangle(G)+ \delta(G))}^2\Big)}}{n\sqrt{\triangle(G) \delta(G)}}. $$
This proves Inequality \ref{E6}.

 The proof of the remaining part of the theorem follows by using similar arguments as in  Theorem \ref{T1}. \qed

Taking $k=1$ in Theorem \ref{T4}, we obtain an upper bound for the signless Laplacian spectral radius as follows.
 \begin{theorem}\label{T5}
Let $G$ be a connected graph with $n$ vertices and $m$ edges. Then
 \begin{equation*}
 q(G)\leq \frac{2m}{n}+\frac{\sqrt{m(n-1)\Big(2\triangle(G)\delta(G)(n-2m)+m{(\triangle(G)+\delta(G))}^2\Big)}}{n\sqrt{\triangle(G)\delta(G)}}
 \end{equation*}
 with equality if and only if $G\cong K_n$.
\end{theorem}

Proceeding and using arguments similar to those used in Theorem \ref{T5}, we get the following lower bound for $L_k(G)$.
 \begin{theorem}\label{T6}
Let $G$ be a connected graph with $n$ vertices and $m$ edges. If $1\leq k\leq n-1$, then
 \begin{equation*}
 L_k(G)\geq \frac{2mk}{n}-\frac{\sqrt{mk(n-k)\Big(2\triangle(G)\delta(G)(n-2m)+m{(\triangle(G)+\delta(G))}^2\Big)}}{n\sqrt{\triangle(G)\delta(G)}}
 \end{equation*}
 with equality if and only if $G\cong K_n$ and $k=n-1$. Equality always holds when $k=n$.
\end{theorem}

\section{Signless Laplacian energy of a graph}

The Laplacian energy of a graph $G$ is defined as $LE(G)=\displaystyle\sum_{i=1}^{n}\Big|\mu_i-\frac{2m}{n}\Big|$. This quantity, which is an extension of graph-energy concept \cite{R7}, has found remarkable chemical applications beyond the molecular orbital theory of conjucated molecules (see \cite{R18}).

In analogy to Laplacian energy, the signless Laplacian energy $QE(G)$ of $G$ is defined as
$$QE(G)=\displaystyle\sum_{i=1}^{n}\Big|q_i-\frac{2m}{n}\Big|.$$
To see the basic properties of this quantity, including various upper and lower bounds, we refer to \cite{R19,R20,R21,R25}.
We start with the following lemma which gives an upper bound for the $Q$-index $q(G)$ of a connected graph $G$ in terms of the order $n$ and size $m$.
\begin{lemma}\label{L9} \emph{ \cite{R5}}
Let $G$ be a connected graph with $n$ vertices and $m$ edges. Then
$$q(G)\leq \frac{2m}{n-1}+n-2$$
with equality if and only if $G$ is $K_{1,n-1}$ or $K_n$.
\end{lemma}

Now, we obtain an upper bound for $QE(G)$ of a connected graph $G$ in terms of the order $n$, size $m$, maximum vertex degree $\triangle(G)$, minimum vertex degree $\delta(G)$ and  $Q$-index $q(G)$ of $G$.
\begin{theorem}\label{T7}
Let $G$ be a connected graph with $n$ vertices and $m$ edges. Then
\begin{equation}\label{E7}
QE(G)\leq \frac{2m}{n(n-1)}+n-2+\sqrt{(n-1)\Big(2m+\frac{n}{4}(\triangle(G)-\delta(G))^2-\Big( q(G)-\frac{2m}{n}\Big)^2\Big)}
\end{equation}
with equality if and only if $G\cong K _n$.
\end{theorem}
\textbf{Proof.} It is easy to see that
$$q_1=q(G)\geq \frac{2m}{n}, ~~~ \displaystyle\sum_{i=1}^{n}\Big|q_i-\frac{2m}{n}\Big|^2=\displaystyle\sum_{i=1}^{n}q^2_i-\frac{4m^2}{n} ~~and~~ \displaystyle\sum_{i=1}^{n}q^2_i=2m+\displaystyle\sum_{i=1}^{n}d^2_i.$$
Using this observations and Lemma \ref{L3}, we get
\begin{align*}
QE(G)&=\displaystyle\sum_{i=1}^{n}\Big|q_i-\frac{2m}{n}\Big|
 = q_1-\frac{2m}{n}+\displaystyle\sum_{i=2}^{n}\Big|q_i-\frac{2m}{n}\Big|\\&
 \leq  q_1-\frac{2m}{n}+\sqrt{(n-1)\displaystyle\sum_{i=2}^{n}\Big|q_i-\frac{2m}{n}\Big|^2}\\&
 = q_1-\frac{2m}{n}+\sqrt{(n-1)\Big(\displaystyle\sum_{i=1}^{n}q^2_i-\frac{4m^2}{n}-\Big( q_1-\frac{2m}{n}\Big)^2\Big)}\\&
 = q_1-\frac{2m}{n}+\sqrt{(n-1)\Big(2m+\displaystyle\sum_{i=1}^{n}d^2_i-\frac{4m^2}{n}-\Big( q_1-\frac{2m}{n}\Big)^2\Big)}\\&
 \leq  q_1-\frac{2m}{n}+\sqrt{(n-1)\Big(2m+\frac{4m^2}{n}+\frac{n}{4}(\triangle(G)-\delta(G))^2-\frac{4m^2}{n}-\Big( q_1-\frac{2m}{n}\Big)^2\Big)}\\&
 \text{(by using Lemma \ref{L5})}\\&
 \leq  \frac{2m}{n-1}+n-2-\frac{2m}{n}+\sqrt{(n-1)\Big(2m+\frac{n}{4}(\triangle(G)-\delta(G))^2-\Big( q_1-\frac{2m}{n}\Big)^2\Big)}\\&
 \text{(by using Lemma \ref{L9})}\\&
 =  \frac{2m}{n(n-1)}+n-2+\sqrt{(n-1)\Big(2m+\frac{n}{4}(\triangle(G)-\delta(G))^2-\Big( q(G)-\frac{2m}{n}\Big)^2\Big)}.
\end{align*}
This proves the required inequality.

Assume that equality holds in Inequality \ref{E7}. Then equality must hold in all the above inequalities, that is, equality must hold simultaneously in Lemmas \ref{L3}, \ref{L5} and \ref{L9}. We consider the following cases.\\
 \textbf{Case 1.} Equality holds in Lemma \ref{L3} if $\Big|q_2-\frac{2m}{n}\Big|= \Big|q_3-\frac{2m}{n}\Big|=\dots=\Big|q_n-\frac{2m}{n}\Big|$.\\
 \textbf{Case 2.} Equality holds in Lemma \ref{L9} if $G$ is either $K_{1,n-1}$ or $K_n$. But $K_{1,n-1}$ does not satisfy Case 1. $K_n$ satisfies Case 1 and also equality holds in Lemma \ref{L5} when $G\cong K_n$ as $K_n$ is a regular graph.\\
 All these arguments show that if equality holds in Inequality \ref{E7}, then $G\cong K_n$.

 Conversely, if $G\cong K_n$, then it is easy to see that the equality holds in Inequality \ref{E7}. \qed

The next lemma due to Cean \cite{R2} gives the upper bound for the sum of the squares of vertex degrees in a graph.

\begin{lemma} \label{L2}\emph{\textbf{\cite{R2}}}
Let $G$ be a graph with $n$ vertices and $m$ edges. Then
$$\displaystyle\sum_{u\in V(G)}d^2_u\leq m\Big(\frac{2m}{n-1}+n-2\Big).  $$
Moreover, if $G$ is connected, then equality holds if and only if $G$ is either a star $K_{1,n-1}$ or a complete graph $K_n$.
\end{lemma}

Proceeding and using arguments similar to Theorem \ref{T7} and using Lemma \ref{L2} in place of Lemma \ref{L5}, we get the following upper bound for $QE(G)$ in terms of order $n$, size $m$ and $Q$-index $q(G)$ of $G$.
  \begin{theorem}\label{T8}
  Let $G$ be a connected graph with $n$ vertices and $m$ edges. Then
\begin{equation*}
QE(G)\leq \frac{2m}{n(n-1)}+n-2+\sqrt{(n-1)\Big(mn+\frac{2m^2(2-n)}{n(n-1)}-\Big( q(G)-\frac{2m}{n}\Big)^2\Big)}
\end{equation*}
with equality if and only if $G\cong K _n$.
  \end{theorem}

\noindent\textbf{Acknowledgement} The research of Prof. S. Pirzada is supported by the SERB-DST research project number CRG/2020/000109. \\

\noindent{\bf Data availibility} Data sharing is not applicable to this article as no datasets were generated or analyzed
during the current study.

\end{document}